\documentclass[showpacs,preprintnumbers,amsmath,amssymb,12pt]{article}
\usepackage{graphicx,amsfonts,amssymb,amsthm,amsmath,graphics,epsfig,pstricks,fancyhdr,fancybox}
\usepackage{dcolumn}
\usepackage{bm}

\newtheorem{thm}{Theorem}[section]

\newtheorem{lem}[thm]{Lemma}
\newtheorem{prop}[thm]{Proposition}

\newcommand{\refeq}[1]{~(\ref{#1})}
\newcommand{\myref}[1]{~\ref{#1}}
\newcommand{\mycite}[1]{~\cite{#1}}
\newcommand{\boxit}[1]{\vbox{\hrule\hbox{\vrule\kern5pt
     \vbox{\kern5pt#1\kern3pt}\kern5pt\vrule}\hrule}}

\newcommand{\sqr}[2]{{\vcenter{\vbox{\hrule height.#2pt
    \hbox{\vrule width.#2pt height#1pt \kern#1pt
    \vrule width.#2pt}\hrule height.#2pt}}}}

\newcommand{\à}{\`a}

\newcommand{\att}[1]{\E\!\left[{#1}\right]}

\newcommand{\pr}[1]{\P\!\left\{{#1}\right\}}

\newcommand{\indim}{\noindent{\bf Proof:}\hspace{0.2cm}}
\newcommand{\findim}{\hfill$\blacksquare$\vspace{0.3cm}\noindent}

\def\RE{{\bm R}}

\def\P{{\bm P}}
\def\Q{{\bf Q}}

\def\E{{\bm E}\,}

\def\OM{\Omega}
\def\om{\omega}

\def\tok{\buildrel k \over \longrightarrow}

\def\va{\emph{rv}}

\def\FD{\emph{cdf}}

\def\ddp{\emph{pdf}}

\def\mmax{k}
\def\Q0{\mathbb{Q}_0}

\begin{document}
\thispagestyle{empty}

\title{\LARGE \textbf{Thou shalt not say \emph{``at random''} in vain:}\\
\LARGE \textbf{Bertrand's paradox exposed}}
\author{\textsc{Nicola Cufaro Petroni}\\
Department of Mathematics and \textsl{TIRES}, University of Bari\\
\textsl{INFN} Sezione di Bari\\
via E. Orabona 4, 70125 Bari, Italy\\
email: \textit{cufaro@ba.infn.it}}

\date{}

\maketitle

\begin{abstract}
\noindent We review the well known Bertrand paradoxes, and we first
maintain that they do not point to any probabilistic inconsistency,
but rather to the risks incurred with a careless use of the locution
\emph{at random}. We claim then that these paradoxes spring up also
in the discussion of the celebrated Buffon's needle problem, and
that they are essentially related to the definition of (geometrical)
probabilities on \emph{uncountably} infinite sets. A few empirical
remarks are finally added to underline the difference between
\emph{passive} and \emph{active} randomness, and the prospects of
any experimental decision
\end{abstract}

{ \small\noindent MSC: 60-01, 60A05, 60E05

\noindent PACS: 02.50.Cw, 01.70.+w

\noindent Key words: Bertrand's paradox; Buffon's needle;
Randomness; Change of measure}

\section{The Bertrand paradoxes}\label{intro}

In the first chapter of his classic treatise\mycite{bert} Joseph
Bertrand dwells for a while on the definition of probability, and in
particular -- in the paragraphs 4-7 -- he remarks that the random
models with an \emph{infinite} number of possible results are prone
to particularly insidious misunderstandings\footnote{``L'infini
n'est pas un nombre; on ne doit pas, sans explication, l'introduire
dans les raisonnements ... Choisir \emph{au hasard} entre un nombre
infini de cas possibles, n'est pas une indication suffisante.''
See\mycite{bert} p.\ 4}. He lists then a few examples of problems
each admitting equally legitimate, but contradictory answers and
suggests then that our questions are \emph{ill posed}, or more
precisely that the required probabilities, based on some \emph{at
random} (\emph{au hasard}) choice, ``sont impossibles \à assigner si
la question n'est pas pr\écis\ée davantage'' (see\mycite{bert} p.\
7). How it will be made clear later, however -- and how it was
likely clear to Bertrand himself -- the crucial point is less the
infinity of the possible outcomes, than their \emph{uncountable}
infinity: a feature shared with other time honored problems, as for
instance that of \emph{Buffon's needle} also discussed later in the
present paper. There is room to argue indeed that even for
\emph{countably} infinite sample spaces the paradoxes do not arise
because the very notion of \emph{at random} -- as long as it is
associated with some idea either of \emph{equiprobability} or of
\emph{uniformity} -- that has not there a straight extension, can be
asymptotically retrieved in a way free from ambiguities. Bertrand on
the other hand, while correctly pointing out that the questions
proposed in his examples are fallacious exactly because our use of
the said locution is too careless, fails to elaborate further on
this point leaving the reader with the odd feeling that something
could be inconsistent in the general notion of randomness. A
negligence extended -- with few notable
exceptions\mycite{gned,mathai} -- also to many of the modern
textbooks that still bother to mention this topic\footnote{See for
instance\mycite{pap} whose final remarks (p.\ 9) are not really
helpful: ``We have thus found not one but three different solutions
for the same problem! One might remark that these solutions
correspond to three different experiments. This is true but not
obvious and, in any case, it demonstrates the ambiguities associated
with the classical definition, and the need for a clear
specification of the outcomes of an experiment and the meaning of
the terms `possible' and `favorable' ''}

The aim of the present paper is then to address this very point:
what are the root and the scope of these seeming inconsistencies?
And in accomplishing our task we will linger first in the
Section\myref{chord} on the example that is widely acknowledged
today as the paradigmatic {Bertrand paradox} because its results
look especially puzzling. We will then proceed in Section\myref{bb}
to extend similar remarks to the Buffon needle problem, and in
Section\myref{count} to argue that while the paradoxes certainly
arise in the event of (geometrical) probabilities defined on
uncountably infinite sets, asymptotically equiprobable countably
infinite sets (as for instance the rational numbers discussed in the
Appendix\myref{randrat}) seem to share the fate of finite sets in
avoiding these ambiguities. In the last Section\myref{conclude} we
will finally conclude by adding a few remarks about the meaning of a
possible experimental discrimination among the different legitimate
solutions

\section{The circle, the triangle and the chord}\label{chord}

\begin{figure}
 \begin{center}
\includegraphics*[width=12cm]{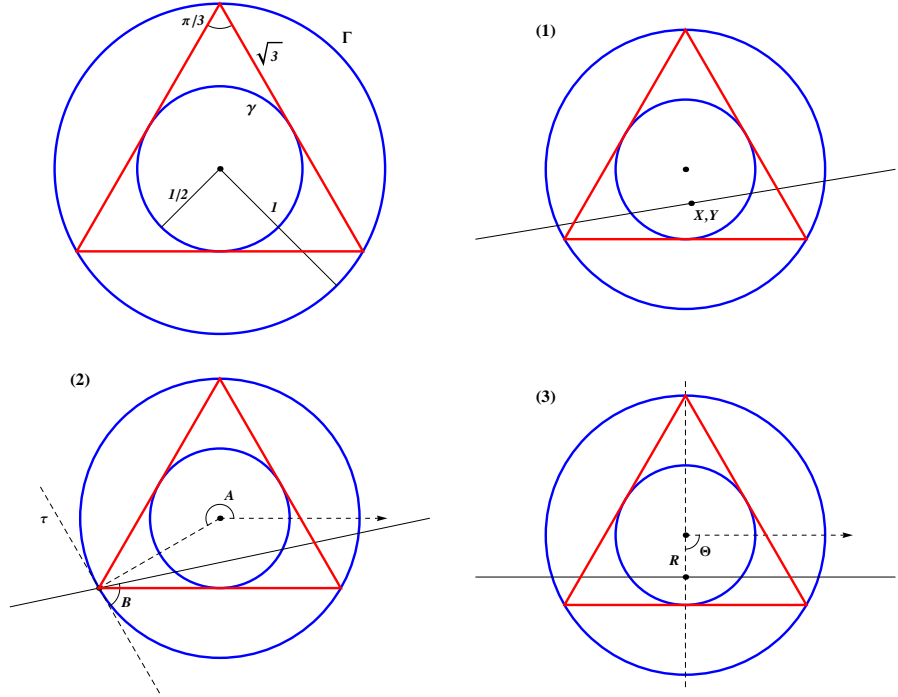}
\caption{Bertrand's paradox}\label{bertrand}
 \end{center}
 \end{figure}

Usually the problem is proposed in the following way: looking at the
Figure\myref{bertrand}, take \emph{at random} a chord on the circle
$\Gamma$ of radius $1$: what is the probability that its length will
exceed that of the edge of an inscribed equilateral triangle (namely
will exceed $\sqrt{3}$)? Three acceptable solutions are possible,
but their answers are all \emph{numerically} different (we always
make reference to the Figure\myref{bertrand}):
\begin{enumerate}
    \item To take a chord at random is equivalent to choose the location of its middle point
    (its orientation would be an aftermath), and to get the chord
    longer than the triangle edge it is necessary and sufficient to
    take this middle point inside the concentric circle
    $\gamma$ with radius $^1/_2$ inscribed in the triangle. The required probability is
    then the ratio between the area $^\pi/_4$ of $\gamma$ and the area $\pi$
    of $\Gamma$, and consequently we have $p_1=\,^1/_4$
    \item By symmetry the position of one chord endpoint along the circle is immaterial to our
    calculations: then, for a given endpoint, the chord length will
    only be contingent on the angle (between $0$ and $\pi$) with the tangent line $\tau$ in the chosen endpoint.
    If then we draw the triangle with one vertex in the chosen endpoint, the chord at random will exceed its edge
    if the angle with the tangent falls between
    $^\pi/_3$ and $^{2\pi}/_3$, and the corresponding probability will be $p_2=\,^1/_3$
    \item Always by symmetry, the random chord direction does not affect the required probability.
    Fix then such a direction, and remark that the chord will exceed $\sqrt{3}$
    if its intersection with the orthogonal diameter falls
    within a distance from the center smaller than $^1/_2$: this happens with probability $p_3=\,^1/_2$
\end{enumerate}
We are then left with three different ($^1/_4,\,^1/_3$ and $^1/_2$),
but equally acceptable answers. To find the paradox origin we must
remember that taking a number \emph{at random} in an uncountably
infinite domain usually means that this number is there
\emph{uniformly} distributed. It is possible to show however (as
also hinted in\mycite{gned}) that what is considered as uniformly
distributed in every single proposed solution can not at the same
time be uniformly distributed in the other two: in other words, in
our three solutions -- by differently choosing \emph{what} is
uniformly distributed -- we surreptitiously adopt three different
probability distributions, and consequently it is not astonishing at
all that the three answers mutually disagree

To be more precise, let us define (see Figure\myref{bertrand}) the
three \va\ (random variable) pairs representing the coordinates
describing the position of our chord in the three proposed
solutions:
\begin{enumerate}
    \item the Cartesian coordinates $(X,Y)$ of the chord middle point
    \item the angles $(A,B)$ respectively giving the position of the fixed endpoint
    and the chord orientation w.r.t.\ the tangent
    \item the polar coordinates $(R,\Theta)$ of the chord-diameter intersection
\end{enumerate}
In every instance however we apparently make the concealed (namely
not explicitly acknowledged) hypothesis that the corresponding pair
of coordinates is uniformly distributed, but these three assumptions
are not mutually consistent, as we will see at once, because they
require three different probability measures on the probability
space where all our \va's are defined. In particular, and by
adopting the notation
\begin{equation*}
    \chi_{[a,b]}(x)=\left\{
                      \begin{array}{ll}
                        1, & \hbox{if $a\le x\le b$;} \\
                        0, & \hbox{else}
                      \end{array}
                    \right.
\end{equation*}
the three solutions respectively assume the following uniform, joint
distributions (see also Figure\myref{densities} for a graphical
account of their respective supports):
 \begin{figure}
 \begin{center}
\includegraphics*[width=16cm]{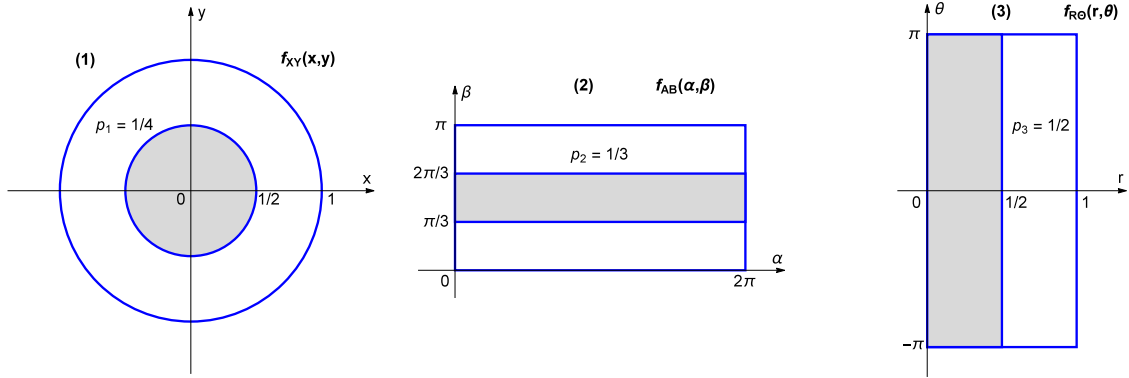}
\caption{Supports of the uniform \ddp's\refeq{XY},\refeq{AB}
and\refeq{RT}: the shaded areas correspond to the three Bertrand
probabilities $p_1, p_2$ and $p_3$}\label{densities}
 \end{center}
 \end{figure}
\begin{enumerate}
    \item the joint, uniform \ddp\ on the unit circle in $\RE^2$
    \begin{equation}\label{XY}
        f_{XY}(x.y)=\frac{1}{\pi}\,\chi_{[0,1]}(x^2+y^2)
    \end{equation}
    for the pair $(X,Y)$: here the two \va's $X$ and $Y$ are \emph{not}
    independent
    \item the joint, uniform \ddp\ on the rectangle $[0,2\pi]\times[0,\pi]$ in $\RE^2$
    \begin{equation}\label{AB}
        f_{AB}(\alpha,\beta)=\frac{1}{2\pi^2}\,\chi_{[0,2\pi]}(\alpha)\chi_{[0,\pi]}(\beta)
    \end{equation}
    for the pair $(A,B)$ with independent components
    \item and finally the joint, uniform \ddp\ on the rectangle $[0,1]\times[-\pi,\pi]$ in $\RE^2$
    \begin{equation}\label{RT}
        f_{R\Theta}(r,\theta)=\frac{1}{2\pi}\,\chi_{[0,1]}(r)\chi_{[-\pi,\pi]}(\theta)
    \end{equation}
     for the pair $(R,\Theta)$ again with independent components
\end{enumerate}
Surely enough, if we would adopt a unique probability space for all
of our three solutions, the three numerical results would be exactly
coincident, but in this case only one of the three \va\ pairs could
possibly be uniformly distributed, while the other joint
distributions should be derived by adopting the well known
procedures established for the functions of \va's (see for
instance\mycite{pap}, Sections $5.2$, $6.2$ and $6.3$). The crucial
point here is that there are in fact some precise transformations
allowing the change from a pair of \va's to the other: by using
these transformations we can show indeed that if a pair is jointly
uniform, then the other two can not have the same property

Without going into the details of every possible combination in our
problem, we will confine ourselves to discuss just the relations
between the solutions (1) and (3). The transformations between the
Cartesian coordinates $(X,Y)$ and the polar ones $(R,\Theta)$ are
well known:
\begin{equation*}
    \left\{
      \begin{array}{l}
        x=r\cos\theta \\
        y=r\sin\theta
      \end{array}
    \right.
    \qquad
    \left\{
      \begin{array}{ll}
        r=\sqrt{x^2+y^2} & \hbox{$\qquad0<r$} \\
        \theta=\arctan\,^y/_x & \hbox{$\,\quad-\pi<\theta\leq\pi$}
      \end{array}
    \right.
\end{equation*}
with a Jacobian determinant
\begin{equation*}
    J(r,\theta)=\left|
                  \begin{array}{cc}
                    ^{\partial r}/_{\partial x} & \,^{\partial r}/_{\partial y} \\
                    ^{\partial\theta}/_{\partial x} & \,^{\partial\theta}/_{\partial y} \\
                  \end{array}
                \right|=
    \left|
       \begin{array}{cc}
         \cos\theta &\; \sin\theta \\
         -\,^1/_r\,\sin\theta &\; ^1/_r\,\cos\theta \\
       \end{array}
     \right|=\frac{1}{r}
\end{equation*}
As a consequence (see for instance\mycite{pap} Section $6.3$), if
$(X,Y)$ have the joint uniform \ddp\refeq{XY}, then the pair
$(R,\Theta)$ will not be uniform and will have instead the \ddp
\begin{equation*}
    f_{R\Theta}^{(1)}(r,\theta)=\frac{r}{\pi}\,\chi_{[0,1]}(r)\chi_{[-\pi,\pi]}(\theta)
\end{equation*}
which is apparently different from the $f_{R\Theta}$ in\refeq{RT}.
By taking advantage of this new distribution $f_{R\Theta}^{(1)}$
(coherent now with the choice of a jointly uniform pair $X,Y$) it is
easy to see that the required probability within the framework of
the solution (3) would be
\begin{equation*}
    p_3^{(1)}=\int_0^\frac{1}{2}\frac{r}{\pi}\,dr\int_{-\pi}^{\pi}d\theta=\frac{1}{4}
\end{equation*}
instead of $p_3=\,^1/_2$, in perfect agreement with the value
$p_1=\,^1/_4$ of the solution (1). Hence the paradox ghosts would
daunt us just as long as we unwittingly suppose that in our three
solutions the coordinates can all be \emph{at once} uniformly
distributed (hiding that under the careless locution \emph{at
random}), and they will disappear instead as soon as we consistently
adopt a unique probability space for all our \va's
 \begin{figure}
 \begin{center}
\includegraphics*[width=9cm]{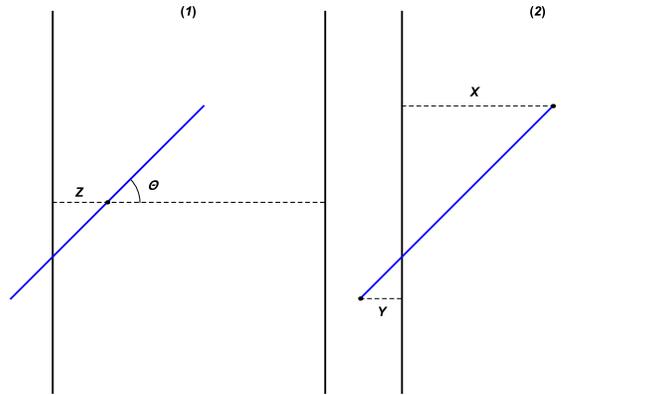}
\caption{Bertrand's paradox for Buffon's needles}\label{buffon}
 \end{center}
 \end{figure}

\section{Bertrand \textit{vs} Buffon}\label{bb}

\begin{figure}
 \begin{center}
\includegraphics*[width=13cm]{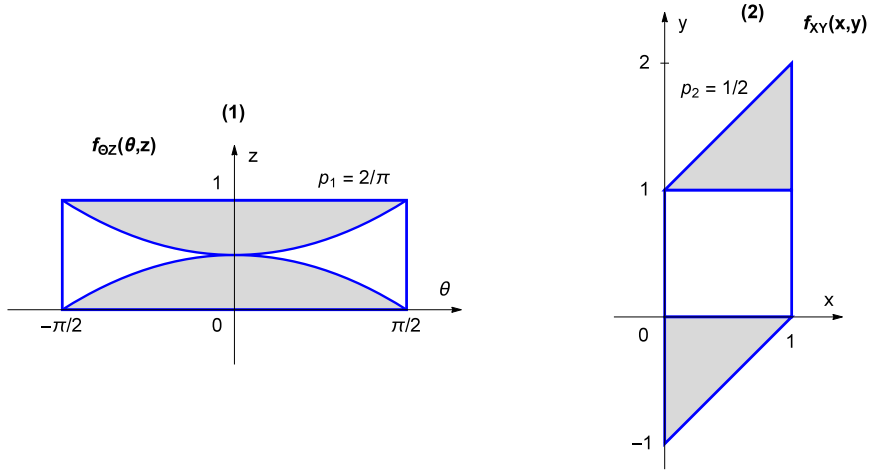}
\caption{Supports of the uniform \ddp's\refeq{TZ} and\refeq{xy}: the
shaded areas correspond to the two Buffon probabilities $p_1$ and
$p_2$}\label{buffondensities}
 \end{center}
 \end{figure}

It is interesting to remark now that, while it is known that by the
turn of the century several different solutions of the Bertrand
question were added\footnote{See for instance\mycite{czu} quoted
in\mycite{mathai}} to the usual three recalled in the previous
section, nobody at our knowledge seem to have noticed that the same
kind of paradoxes does in fact appear also in the discussion of the
celebrated Buffon needle problem. In its simplest
version\footnote{For a more complete discussion see for
instance\mycite{gned} and\mycite{shir}} a needle of unit length is
thrown \emph{at random} on a table where a few parallel lines are
drawn at a unit distance: what is the probability that the needle
will lie across one of these lines? In the classical answer to this
question, since the lines are drawn periodically on the table, it
will be enough to study the problem with only two lines by supposing
that the needle center does fall between them. The position of the
said center along the direction of the parallel lines is also
immaterial. The needle position is then defined by just two
\emph{\va}'s: the distance $Z$ of its center from the left line, and
the angle $\Theta$ between the needle and a perpendicular to the
parallel lines (see $(1)$ in Figure\myref{buffon}). That the needle
is thrown \emph{at random} here means that the pair of \emph{\va}'s
$\Theta,Z$ is uniform in $[-{\pi\over2},{\pi\over2}]\times[0,1]$,
namely that their joint \ddp\ is
\begin{equation}\label{TZ}
     f_{\Theta Z}(\theta,z)=\frac{1}{\pi}\,\chi_{[-{\pi\over2},{\pi\over2}]}(\theta)\,\chi_{[0,1]}(z)
\end{equation}
while, with $-{\pi\over2}\leq \theta\leq{\pi\over2}$, the needle
will lie across a line either when $x\leq\,^1/_2\cos \theta$, or
when $x\geq1-\,^1/_2\cos \theta$ (see again $(1)$ in
Figure\myref{buffon}). Just by inspecting the \ddp\ $(1)$ in
Figure\myref{buffondensities} it is then easy to find out that the
required probability is
\begin{equation*}
    p_1={2\over\pi}\int_{-\pi/2}^{\pi/2}\frac{\cos \theta}{2}\,d\theta={2\over\pi}
\end{equation*}

\noindent In the spirit of the Bertrand paradoxes, however, we can
give a different answer to the Buffon question (see $(2)$ in
Figure\myref{buffon}): the needle position is now identified by
looking first to its (vertically) upper end, and by recording its
distance $X$ from the left line. Then we consider where its other
(lower) end falls and we mark down its distance $Y$ from the same
left line. If the said left line is in the origin of a horizontal
axis, it is apparent that for every value $0\le x\le1$ of $X$, the
possible values of $Y$ will be between $x-1$ and $x+1$ (because
apparently $|x-y|\le1$), and in this framework to throw the needle
\emph{at random} will mean that the joint distribution of $X,Y$ is
uniform in the domain shown in $(2)$ of the
Figure\myref{buffondensities}, namely
\begin{equation}\label{xy}
     f_{XY}(x,y)=\frac{1}{2}\,\chi_{[0,1]}(x)\,\chi_{[0,1]}(|x-y|)
\end{equation}
On the other hand it is apparent that for every $0\le x\le1$ the
needle will cross a line when either $x-1\le y\le0$, or $0\le y\le
x+1$, so that the required probability will correspond to the shaded
area in $(2)$ of Figure\myref{buffondensities} and hence now
$p_2=\,^1/_2$. These remarks show then that also the Buffon needle
problem is not impervious to paradoxes, and this could be more than
just a trifle because of its peculiar \emph{experimental} status, as
will be argued in the subsequent Section\myref{conclude}

\section{Counting and measuring}\label{count}

Since the Bertrand paradoxes arise from a careless use of the
locution \emph{taking at random}, it would expedient to recall once
more that this is not understood here as the drawing of some outcome
$\om$ out of a set $\OM$ according to some arbitrary probability
measure, but stands rather for assuming that there is no reason to
think that there are preferred outcomes $\om\in\OM$, these being
supposed instead to be equally likely. For sets of numbers this kind
of {randomness} is enforced either by sheer \emph{equiprobability}
(on the finite sets, by \emph{counting}), or by distribution
\emph{uniformity} (on the bounded, Lebesgue measurable, uncountable
sets, by \emph{measuring}). On the other hand infinite, countable
sets and unbounded, uncountable sets are both excluded from these
egalitarian probability attributions because their elements can be
made neither equiprobable (with a non vanishing probability), nor
uniformly distributed (with a non vanishing probability density).
This in particular seems to point to the fact that infinite,
countable sets should be barred even from discussing the Bertrand
problems because, for instance, we can not make all the natural
numbers equiprobable with a non zero probability (neither by
counting, nor by measuring), while the paradoxes are essentially
based on a misunderstanding about this kind of randomness. In these
cases however it is possible to start with some proper (not
equiprobable) probability distribution, and then to make them ever
closer -- in a suitable, approximate sense -- to an equiprobable
one: we will then speak of \emph{asymptotic equiprobability} (see
Appendix\myref{randrat} and\mycite{cuf} for more details). This
notion will be used here to argue that the Bertrand paradoxes arises
exclusively in connection to \emph{uncountable} infinite sets of
numbers, but neither for finite, nor for countably infinite ones

To clarify this last point it would be expedient to consider
another, more simple case among the Bertrand examples
(see\mycite{bert} p.\ 4): if we ask what is the probability that a
\emph{real number} $x$ chosen at random between $0$ and $100$ is
larger than $50$, our natural answer is $^1/_2$. Since however the
real numbers between $0$ and $100$ are also bijectively associated
to their squares between $0$ and $10\,000$, we also instinctively
feel that our question should be equivalent to ask for the
probability that the square of a real random number turns out to be
larger than $50^2=2\,500$. When however we think of taking this last
number \emph{at random} between $0$ and $10\,000$, instinctively
again we are inclined to answer that the probability of exceeding
$2\,500$ should now be $^3/_4$ instead of $^1/_2$. The two problems
look equivalent, but their two intuitive answers (apparently both
legitimate) are different

Predictably the paradox resolution is similar to that of
Section\myref{chord}:  we would readily concede that the probability
to exceed $50$ for a real number $X$ taken \emph{at random} in
$[0,100]$ is $p_1=\,^1/_2$. When however we ask for the probability
that $X^2$ taken \emph{at random} in $[0,10\,000]$ exceeds
$50^2=2\,500$, we surreptitiously change our measure by supposing
that now $X^2$ is uniform in $[0,10\,000]$ and we find
$p_2=\,^3/_4$. But the fact is -- as in the previous example -- that
if $X$ is uniform in $[0,100]$, then $X^2$ can not be uniform in
$[0,10\,000]$, and vice-versa. More precisely, if the \ddp\ of $X$
is the uniform
\begin{equation*}
    f_X(x)=\frac{1}{100}\,\left\{
             \begin{array}{cl}
               1 & \quad\hbox{for $0\le x\le100$,} \\
               0 & \quad\hbox{else}
             \end{array}
           \right.
\end{equation*}
then the corresponding, non uniform \ddp\ of $Y=X^2$ is (see
again\mycite{pap} Section $5.2$)
\begin{equation*}
    f_Y(y)=\frac{1}{200\sqrt{y}}\,\left\{
             \begin{array}{cl}
               1 & \quad\hbox{for $0\le y\le10\,000$,} \\
               0 & \quad\hbox{else}
             \end{array}
           \right.
\end{equation*}
and of course the paradox disappears because now, in agreement with
$p_1=\,^1/_2$, we would have
\begin{equation*}
    p_2^{(1)}=\int_0^{2\,500}f_Y(y)\,dy=\frac{1}{200}\int_0^{2\,500}\frac{dy}{\sqrt{y}}=\frac{1}{2}
\end{equation*}

It is easy to see moreover that the paradox does not show up at all
when we consider the \emph{finite} version of this problem: if we
ask for the probability ($p_1=\,^1/_2$) of choosing \emph{at random}
an \emph{integer number} $n$ larger than $50$ among the
(\emph{equiprobable}) integer numbers from $1$ to $100$, we would in
fact recover the same answer ($p_2=\,^1/_2$) also by asking to
calculate the probability of choosing \emph{at random} a number
larger than $2\,500$ among the squared integers
$1,4,9,\ldots,10\,000$, because now our set is again constituted of
just $100$ equiprobable elements and there are $50$ larger than
$2\,500$. The crucial difference with the previous continuous
version of the problem is that in the case of finitely many
(equiprobable) possible results we just \emph{enumerate} the
\emph{favorable} and the \emph{possible} items (a situation not
changed by squaring the numbers), while for the continuous real
numbers (geometric probabilities) we \emph{compare} the length of
the intervals: all is contingent indeed on the difference between
\emph{counting} and \emph{measuring}

This situation, albeit trickier, is not essentially changed for
\emph{countably infinite} possible outcomes, but for the fact that
in this case they can not be made strictly equiprobable by direct
enumeration. Take for instance the problem of asking for the
probability of choosing \emph{at random} a \emph{rational number}
$q=\,\!^n/_m$ larger than $50$ among the rational numbers in
$[0,100]$ that are famously an infinite, countable set everywhere
dense among the real numbers. It is shown in the
Appendix\myref{randrat} that -- in agreement with our intuition --
this probability tends to $p_1=\,^1/_2$ when, going around the
problem of actually enumerating them, the rational numbers in
$[0,100]$ are made asymptotically equiprobable. If however we
subsequently ask to calculate the probability of choosing \emph{at
random} a number larger than $2\,500$ among the squared rational
numbers $q^2=\,\!^{n^2}\!/_{m^2}$ in $[0,10\,000]$ we in fact
recover the same answer $p_2=\,^1/_2$ with no possible ambiguity
because now -- in a way recalling the case of integer numbers -- we
must make asymptotically equiprobable not all the rational numbers
in $[0,10\,000]$, but rather only those that are squares of rational
numbers. Not every rational $q$ is indeed a squared rational, and
the (tiny, but non zero) probability of $q^2=\,\!^{n^2}\!/_{m^2}$
being larger than $2\,500$ exactly coincides with that of
$q=\,\!^n/_m$ being larger than $50$: there is no possible sharing
of probability with the infinitely many other (non squared)
rationals that in any case would never show up in the process of
drawing at random squared rationals

The critical feature -- common to both the finite and the countably
infinite cases -- appears here to be the possibility that every
single element be endowed with its own individual \emph{non
vanishing} probability that it also carries with him in every
conceivable one-to-one transformation: no amount of probability must
indeed be shared with numerical results other than the transformed
results, and the distribution remain the same in the transformed
sample space. A situation totally at variance with that of a
\emph{geometrical} probability on a uncountably infinite set where
every sample is usually entitled only to strictly zero probability
and the transformations usually imply a stretching of the
probability measure

\section{An almost empiric conclusion}\label{conclude}

Going back to the example discussed in the Section\myref{chord}, we
have argued that the paradoxes ``disappear as soon as we
consistently adopt a unique probability space for all our \va's'',
so that -- paradoxes notwithstanding -- there are no possible formal
inconsistencies within our overall probabilistic framework. But this
is sheer mathematics, and we are left anyway with three (or more)
possible, coherent and perfectly legitimate, probabilistic models
giving rise to three numerically different results: which one is
\emph{true}, in the sense that it corresponds to the \emph{physical
reality}? This problem of course can not be solved with a
calculation, and should instead be settled -- if possible -- by
comparing the solutions proposed with some empiric result. In other
words one should simply \emph{perform the experiment of choosing at
random a chord on a circle} (or something equivalent) in order to
compare then its statistics with the calculations: something that to
date, at our knowledge, has not yet been done once and for all, and
it is not even exactly tackled in recent contributions\mycite{aerts}
where, for all the emphasis on solving this {hard part} of the
paradox, the discussion seems to be essentially restricted to the
mathematical models

In this vein we will add here just a few final remarks: first, the
possible experiment can not definitely be a \emph{simulated} one
performed on a computer. In this case indeed our experimenter should
\emph{a priori} choose one of the three models to program his
computer to produce a particular pair of uniformly distributed
coordinates: but in so doing he would have already decided the
outcomes of the experiment that coherently will now confirm the
chosen model: this apparently will prove nothing. Second, it is
possible that even in some real, physical experiment the outcome can
be influenced by the choice of what exactly we decide to
measure\mycite{aerts}: different experimental settings could point
to different facets of the physical reality, and after all the
probability is not listed among the concrete things of this world
(see for instance the unconventional viewpoint displayed
in\mycite{defin}, vol.\ 1, \emph{Preface}) representing rather the
state of our information. Third, the previous remark also lays bare
the difference between experiments where we study a \emph{passive
randomness} produced by an independent external world (think for
instance to the statistics of the usual empirical measures, or even
to quantum mechanics: randomness is there completely outside of our
control and show characters that must be discovered rather than
produced), in contrast with an \emph{active randomness} where we try
to produce some kind of previously planned events, namely to
empirically enforce some idea of randomness as in both the Bertrand
paradox and the Buffon needle, and from a different standpoint in
every computer simulation. These two sorts of randomness appear
epistemologically rather different and this dissimilarity could be
well worth of further inquiry

We can not refrain however from pointing out in the end that the
presence of Bertrand-type paradoxes even in the discussion of
Buffon's needle sheds a different light on this problem: it is known
indeed that the classical calculation of the Buffon needle, giving
the probability $p_1=\,^2/_\pi$, also stimulated several experiments
used to get an empirical determination of the approximate value of
the number $\pi$ (four of these tests dated from 1850 to 1901 are
listed for instance in\mycite{gned}), and famously known as
pioneering examples of the \emph{Monte Carlo} method. Despite a few
reservations about the reliability of these results\mycite{gned}, it
is striking that all these four experiments point to a number in the
range between $3.14$ and $3.16$, while our second, proposed
alternative solution with $p_2=\,^1/_2$ would require results
clustering around $4.00$. It is possible -- as suggested above --
that the quoted results are biased by some unaware bent toward
$\pi$, but if confirmed they would suggest that there could be an
empirical meaning in the locution \emph{at random} because, at least
in the case of Buffon's needle, the experiments appear to be able to
favor one among several formally legitimate solutions. But it is
also apparent by now that a possible answer to these questions would
lie outside the reach of this paper, so that for the time being we
will stay content with having just clarified the meaning and the
scope of the Bertrand paradoxes and their link to the Buffon needle,
by leaving to future inquiries the practical task of empirically
deciding among the mathematical models

\vspace{10pt}

\noindent \textbf{Acknowledgements}: The author would like to thank
A.\ Andrisani, S. Pascazio and C.\ Sempi for stimulating discussions
and suggestions

\begin{appendix}

\section{Taking rational numbers at random}\label{randrat}

Rational numbers are famously countable, and hence they can be put
in a sequence. Since however they are a dense subset of the real
numbers, every rational number is a cluster point, and hence no
sequence encompassing all of them can ever converge, not to say be
monotone. In any case their countability certainly allows the
allotment of discrete distributions with non vanishing probabilities
for every item: since they are infinite, however, they can never be
exactly \emph{equiprobable}. In this appendix (further details are
available in\mycite{cuf}) we will outline a procedure to give
distributions on the rationals in $[0,1]$, a set that we will
shortly denote as $\Q0=\mathbb{Q}\cap[0,1]$, and we will investigate
if and how they can be made \emph{asymptotically equiprobable}

It is however advisable to assert right away that the distribution
of a \va\ $Q$ taking values in $\Q0$ must anyhow be of a discrete
type, allotting (possibly non vanishing) probabilities to the
individual rational numbers $q\in\Q0$: conceivable continuous set
functions -- namely with continuous, albeit perhaps not absolutely
continuous, \FD\ (cumulative distribution function) -- would turn
out to be not countably additive, and hence would not qualify as
measures, not to say as probability distributions. Every continuous
\FD\ for $Q$ would indeed entail that at the same time
$\pr{Q=q}=0,\;\forall q\in\Q0$, and $\pr{Q\in\Q0}=1$, while $\Q0$
apparently is the countable union of the disjoint, negligible sets
$\{q\}$: in plain conflict with the countable additivity. This in
particular rules out for the numbers in $\Q0$ the possibility of
being \emph{uniformly distributed} (an imaginable surrogate of
equiprobability evoked by the density of the rational numbers): this
property would in fact require a continuous \FD\

Taking advantage now of the well known diagram used to prove the
countability of the rational numbers, we will consider two dependent
\va's $M$ and $N$ with integer values
\begin{equation*}
    m=1,2,\ldots \qquad\qquad n=0,1,2,\ldots,m
\end{equation*}
and acting respectively as denominator and numerator of the random
rational number $Q=\,^N/_M\in[0,1]$. As a consequence $Q$ will take
the values $q=\,^n/_m$ arrayed in a triangular scheme as in
Table\myref{rationals}.
\begin{table}
\begin{center}
    \begin{tabular}{|c|cccccccccc}
  &  &  &  &  & $\bm n$ &  &  &  & & \\
  \hline
    $\bm m$ & $0$ & $1$ & $2$ & $3$ & $4$ & $5$ &$6$  & $7$ & $8$ & $\ldots$ \\
  \hline
  % after \\: \hline or \cline{col1-col2} \cline{col3-col4} ...
  $1$  & $0$ & $1$ &  &  &  &  &  & & &\\
  $2$  & $0$ & $^1/_2$ & $1$ &  &  &  &  & & &\\
  $3$  & $0$ & $^1/_3$ & $^2/_3$ & $1$ &  &  &  & & &\\
  $4$  & $0$ & $^1/_4$ & $^2/_4$ & $^3/_4$ & $1$ &  &  & & &\\
  $5$  & $0$ & $^1/_5$ & $^2/_5$ & $^3/_5$ & $^4/_5$ & $1$ &  & & & \\
  $6$  & $0$ & $^1/_6$ & $^2/_6$ & $^3/_6$ & $^4/_6$ & $^5/_6$ & $1$ & & &\\
  $7$  & $0$ & $^1/_7$ & $^2/_7$ & $^3/_7$ & $^4/_7$ & $^5/_7$ & $^6/_7$ & $1$ & &\\
  $8$  & $0$ & $^1/_8$ & $^2/_8$ & $^3/_8$ & $^4/_8$ & $^5/_8$ & $^6/_8$ & $^7/_8$ & $1$ &\\
  $\vdots$  & $\vdots$ &  &  &  &  &  &  & & & $\ddots$\\
  %\hline
\end{tabular}
\caption{Table of rational numbers $q=\,\!^n/_m$ with repetitions:
many fractions are reducible to canonical forms already present in
earlier positions}\label{rationals}
\end{center}
\end{table}
It is apparent however that in this way every rational number $q$
shows up infinitely many times due to the presence of reducible
fractions: for instance -- with the usual notation for repeating
decimals -- we have
\begin{equation*}
    0.5=\,\!^1/_2=\,\!^2/_4=\,\!^3/_6=\ldots\qquad 0.\overline{3}=\,\!^1/_3=\,\!^2/_6=\ldots\qquad 0.75=\,\!^3/_4=\,\!^6/_8=\ldots
\end{equation*}
If we adopt however the notation
\begin{equation*}
    q\doteq\,\!^n/_m
\end{equation*}
to indicate that $^n/_m$ is the unique irreducible representation of
a rational number $q$, namely that $n$ and $m$ are co-primes, the
previous examples will be listed as
\begin{equation*}
    0.5\doteq\,\!^1/_2\qquad\quad0.\overline{3}\doteq\,\!^1/_3\qquad\quad0.75\doteq\,\!^3/_4
\end{equation*}
By introducing now a joint distributions of $N$ and $M$
\begin{align*}
  &  \pr{M=m} \qquad\qquad\qquad\; m=1,2,\ldots \\
  &  \pr{N=n\left|M=m\right.} \qquad\quad \,n=0,1,2,\ldots,m \\
  &  \pr{N=n,M=m}=\pr{N=n\left|M=m\right.}\,\pr{M=m}
\end{align*}
we will have for every rational $0\le q\doteq\,\!^n/_m\le1$ the
discrete distribution
\begin{eqnarray}\label{discrD}
    \pr{Q=q}&=&\sum_{\ell=1}^\infty\pr{N=\ell n,M=\ell m}\nonumber\\
  &=&\sum_{\ell=1}^\infty\pr{N=\ell n\left|M=\ell m\right.}\,\pr{M=\ell m}
\end{eqnarray}
This also allows to define the \FD\ of $Q$ as (here of course
$x\in\RE$)
\begin{eqnarray}
  F_Q(x) &=& \pr{Q\le x}=\pr{N\le Mx}=\sum_{m=1}^\infty\pr{N\le
  mx\left|M=m\right.}\pr{M=m}\nonumber\\
   &=&\sum_{m=1}^\infty F_N(mx|M=m)\pr{M=m}\label{cdf1}
\end{eqnarray}
and hence also the probability of $Q$ falling in $(a,b]$ for $0\le
a<b\le1$ real numbers:
\begin{eqnarray}
  \pr{a< Q\le b}&=&F_Q(b)-F_Q(a)\nonumber\\
  &=&\sum_{m=1}^\infty \big[F_N(mb|M=m)-F_N(ma|M=m)\big]\pr{M=m}\label{probint1}
\end{eqnarray}
Notice that the conditional \FD\ of $N$ can also be given as
\begin{eqnarray}
  F_N(x|M=m)&=&\pr{N\le
  x\left|M=m\right.}=\sum_{n=0}^m\pr{N=n\left|M=m\right.}\vartheta(x-n)\nonumber\\
   &=&\sum_{n=0}^{\lfloor
   x\rfloor}\pr{N=n\left|M=m\right.}\label{condcdf}
\end{eqnarray}
where
\begin{equation*}
    \vartheta(x)=\left\{
                   \begin{array}{ll}
                     1 & \qquad x\ge0 \\
                     0 & \qquad x<0
                   \end{array}
                 \right.
\end{equation*}
is the Heaviside function, while for every real number $x$, the
symbol $\lfloor x\rfloor$ denotes the \emph{floor} of $x$, namely
the greatest integer less than or equal to $x$. As a consequence the
equations\refeq{cdf1} and\refeq{probint1} also take the
form\begin{eqnarray}
  F_Q(x)\!\! &=&\!\! \sum_{m=1}^\infty \pr{M=m}\sum_{n=0}^{\lfloor
   mx\rfloor}\pr{N=n\left|M=m\right.} \label{cdf}\\
  \!\!\!\!\!\!\!\!\pr{a< Q\le b}\!\! &=&\!\! \sum_{m=1}^\infty\pr{M=m}(1-\delta_{\lfloor ma\rfloor,\lfloor mb\rfloor})\!\!\! \sum_{n=\lfloor
   ma\rfloor+1}^{\lfloor
   mb\rfloor}\!\!\!\!\pr{N=n\left|M=m\right.} \label{probint}
\end{eqnarray}
where the Kronecker delta takes into account the fact that when
$\lfloor mb\rfloor=\lfloor ma\rfloor$ the difference vanishes, so
that $\lfloor mb\rfloor\ge\lfloor ma\rfloor+1$.

Let us suppose now for simplicity that for a given denominator
$m\ge1$ the $m+1$ possible values of the numerator $n=0,1,\ldots,m$
are equiprobable in the sense that
\begin{equation*}
  \pr{N=n\left|M=m\right.}=\frac{1}{m+1}\qquad\quad n=0,1,\ldots,m
\end{equation*}
In this case for the distribution, with $n,m$ co-primes and $0\le
n\le m$, from\refeq{discrD} we have
\begin{equation}\label{ratprob}
  \pr{Q=q} =
  \sum_{\ell=1}^\infty\frac{\pr{M=\ell m}}{\ell m+1}\qquad\quad q\doteq\,\!^n/_m
\end{equation}
while for the \FD\refeq{cdf} we have from\refeq{condcdf}
\begin{eqnarray}
  F_N\left(mx\left|M=m\right.\right)&=&
  \frac{1}{m+1}\sum_{n=0}^m\vartheta\left(mx-n\right)=\left\{
                                                                 \begin{array}{ll}
                                                                   0 &\quad \hbox{$x<0$} \\
                                                                   \frac{\lfloor mx\rfloor+1}{m+1} &\quad \hbox{$0\le x<1$} \\
                                                                   1 &\quad \hbox{$1\le x$}
                                                                 \end{array}
                                                               \right.\nonumber
  \\
  F_Q(x) &=&\sum_{m=1}^\infty\frac{\pr{M=m}}{m+1}\sum_{n=0}^m\vartheta\left(mx-n\right)\nonumber\\
                                                               &=&\left\{
                                                                 \begin{array}{ll}
                                                                   0 &\quad \hbox{$x<0$} \\
                                                                   \sum_{m\ge1}\pr{M=m}\frac{\lfloor mx\rfloor+1}{m+1} &\quad \hbox{$0\le x<1$} \\
                                                                   1 &\quad \hbox{$1\le x$}
                                                                 \end{array}
                                                               \right.\label{ratcdf}
\end{eqnarray}
and the probability\refeq{probint} with $0\le a<b\le1$ becomes
\begin{equation}\label{ratint}
  \pr{a< Q\le b}=\sum_{m=1}^\infty\pr{M=m}\frac{\lfloor mb\rfloor-\lfloor ma\rfloor}{m+1}
\end{equation}
By denoting now as $p_m=\pr{M=m}$ the distribution of $M$, and as
$s=\sup_m p_m$ the supremum of all its values, let us take now a
sequence of denominators $\{M_k\}_{k\ge1}$ with distributions
$\{p_m(k)\}_{k\ge1}$, and with $s_k$ vanishing for $k\to\infty$ in
such a way that
\begin{equation}\label{log}
    \lim_k\, s_k\ln k=0
\end{equation}
In other words we consider a sequence of distributions that are
increasingly (and uniformly) flattened toward zero, so that the
denominators too are increasingly equiprobable. Ready examples of
these sequences with $k=1,2,\ldots$ are for instance the
\emph{geometric} distributions
\begin{equation*}
    p_m(k)=w_k(1-w_k)^{m-1}\quad\qquad m=1,2,\ldots
\end{equation*}
with infinitesimal $w_k$, and the \emph{Poisson} distributions
\begin{equation*}
    p_m(k)=e^{-\lambda_k}\frac{\lambda_k^{m-1}}{(m-1)!}\quad\qquad m=1,2,\ldots
\end{equation*}
with divergent $\lambda_k$
 \begin{lem}
Within the previous notations and conditions we have
\begin{equation}\label{lem}
   \mu_k=\att{\,^1/_{M_k}}=\sum_{m=1}^{\infty}\frac{p_m(k)}{m}\tok0
\end{equation}
 \end{lem}
 \indim
The positive series defining $\mu_k$ is certainly convergent because
\begin{equation*}
    \mu_k=\sum_{m=1}^{\infty}\frac{p_m(k)}{m}<\sum_{m=1}^{\infty}p_m(k)=1
\end{equation*}
and hence we can always write
\begin{equation*}
    \mu_k=\sum_{m=1}^{\infty}\frac{p_m(k)}{m}=\sum_{m=1}^{k}\frac{p_m(k)}{m}+R_k
\end{equation*}
where
\begin{equation*}
    R_k=\sum_{m=k+1}^{\infty}\frac{p_m(k)}{m}\tok0
\end{equation*}
is an infinitesimal remainder. On the other hand, under our stated
conditions
\begin{equation*}
    \sum_{m=1}^{k}\frac{p_m(k)}{m}<s_k\sum_{m=1}^{k}\frac{1}{m}=s_k
    H_k
\end{equation*}
where $H_k$ denotes the $k^{th}$ \emph{harmonic number}, namely the
sum of the reciprocal integers up to $^1/_k$: it is well\ known
(\cite{grad} $\bm{0.131}$) that for $k\to\infty$ the $H_k$ grow as
$\ln k$, so that from\refeq{log} we have $s_kH_k\tok0$, and finally
$\mu_k=s_kH_k+R_k\tok0$
 \findim
 \begin{prop}
If $Q=\,\!^N/_M$ and $F_Q(x)$ is its \emph{\FD}, then, within the
notation and conditions outlined above, we have
\begin{align}
  & \lim_{\mmax}\pr{Q=q}=0\qquad\qquad\lim_{\mmax}\pr{a<Q\le
               b}=b-a\label{equiprob}\\
%\end{equation}
%and in particular
%\begin{equation}
   & \qquad\qquad\quad\lim_{\mmax}F_Q(x)=\left\{\begin{array}{ll}
                                       0 &\quad \hbox{$x<0$} \\
                                       x &\quad \hbox{$0\le x<1$} \\
                                       1 &\quad \hbox{$1\le x$}
                              \end{array}
                              \right.\label{unifcdf}
\end{align}
 \end{prop}
 \indim
Since our series have positive terms the first result
in\refeq{equiprob} follows from\refeq{ratprob} and\refeq{lem}
because, with $q\doteq\,\!^n/_j$
\begin{equation*}
    \pr{Q=q}=\sum_{\ell=1}^\infty\frac{p_{\ell j}(k)}{\ell j+1}<\sum_{m=1}^\infty\frac{p_{m}(k)}{m+1}
    <\sum_{m=1}^\infty\frac{p_{m}(k)}{m}=\mu_k\tok0
\end{equation*}
As for the second result in\refeq{equiprob}, since for every real
number $x$ it is $x-1\le\lfloor x\rfloor\le x$, for every
$\mmax=1,2,\ldots$, and $0\le a<b\le1$, we have from\refeq{ratint}
\begin{equation*}
   \sum_{m=1}^{\infty}p_m(k)\frac{m(b-a)-1}{m+1} \le\,\pr{a<Q\le b}\,\le\sum_{m=1}^{\infty}p_m(k)\frac{m(b-a)+1}{m+1}
\end{equation*}
namely
\begin{equation*}
   b-a+(a-b-1)\sum_{m=1}^{\infty}\frac{p_m(k)}{m+1}\le\,\pr{a<Q\le b}\,\le b-a+(a-b+1)\sum_{m=1}^{\infty}\frac{p_m(k)}{m+1}
\end{equation*}
so that, since $a-b-1\le0$ and $a-b+1\ge0$, it is
\begin{equation*}
    b-a+(a-b-1)\mu_k\le\,\pr{a<Q\le b}\,\le b-a+(a-b+1)\mu_k
\end{equation*}
The second result\refeq{equiprob} follows then from\refeq{lem}. In a
similar way we finally find for\refeq{unifcdf} that
\begin{equation*}
   \sum_{m=1}^{\infty}p_m(k)\frac{mx}{m+1} \le F_Q(x)\le\sum_{m=1}^{\infty}p_m(k)\frac{mx+1}{m+1}\qquad\quad 0\le x\le1
\end{equation*}
namely
\begin{equation*}
   x-x\sum_{m=1}^{\infty}\frac{p_m(k)}{m+1}\le F_Q(x)\le x+(1-x)\sum_{m=1}^{\infty}\frac{p_m(k)}{m+1}
\end{equation*}
and hence
\begin{equation*}
   x-x\,\mu_k< F_Q(x)< x+(1-x)\mu_k
\end{equation*}
so that the result again follows from\refeq{lem}
 \findim

\noindent From this proposition we see that in the limit
$\mmax\to\infty$, while the probability of every single rational
number rightly vanishes, the probability of these numbers lumped
together in intervals does not: a behavior highly reminiscent of
what happens to continuously distributed \emph{real} \va's. For the
reasons presented at the beginning of this appendix, however, the
previous result by no means imply that we can implement a uniform
limit distribution on $\Q0$ (as we said: there is not such a thing),
but it rather suggests that our random rational numbers $Q$ -- at
least for denominators $m$ distributed in a fairly flat way, and
numerators $n$ conditionally equiprobable between $0$ and $m$ --
asymptotically behave as uniformly distributed in $[0,1]$, and hence
they quite reasonably correspond to our intuitive idea of
\emph{taking rational numbers at random.}

\end{appendix}

\end{document}